\documentclass{birkjour}
\usepackage{amssymb}
\usepackage{color}
\usepackage{hyperref}
\hypersetup{
    linktocpage=true,
    colorlinks=true,
    linkcolor=blue,
    filecolor=magenta,      
    urlcolor=cyan,
    pdftitle={Extreme point methods in the study of isometries on certain non-commutative spaces},
    bookmarks=true,
   citecolor=red    }
 \newtheorem{thm}{Theorem}[section]
 \newtheorem{cor}[thm]{Corollary}
 \newtheorem{lem}[thm]{Lemma}
 \newtheorem{prop}[thm]{Proposition}
 \theoremstyle{definition}
 
 \theoremstyle{remark}
 \newtheorem{rem}[thm]{Remark}
 
 \numberwithin{equation}{section}

\newcommand{\tma}{S(\mathcal{A},\tau)} 
\newcommand{\nmb}{S(\mathcal{B},\nu)} 

\newcommand{\A}{\mathcal{A}}
\newcommand{\B}{\mathcal{B}}

\newcommand{\pa}{\mathcal{P}(\mathcal{A})}    
\newcommand{\pb}{\mathcal{P}(\mathcal{B})}    
\newcommand{\paf}{\mathcal{P}(\mathcal{A})_f}    
\newcommand{\pbf}{\mathcal{P}(\mathcal{B})^f}    

\newcommand{\dfs}[2]{d\left({#1}\right)({#2})}   

\newcommand{\svf}[1]{\mu_{{#1}}}  

\newcommand{\svft}[2]{\mu_{{#1}}({#2})}  

\newcommand{\norm}[1]{\bigl\Vert{#1} \bigr\Vert}

\newcommand{\norminf}[1]{\bigl\Vert{#1} \bigr\Vert_{\infty}}

\newcommand{\ra}{\rightarrow}

\newcommand{\id}{\mathbf{1}} 

\newcommand{\seq}[1]{({#1}_n)_{n=1} ^{\infty}}

\newcommand{\net}[1]{\{{#1}_{\lambda}\}_{\lambda \in \Lambda}}
   
\newcommand{\sumx}[2]{\underset{#1}{\overset{#2}{\sum}}}

\newcommand{\rax}[1]{\overset{#1}{\rightarrow}}
\newcommand{\raxl}[1]{\overset{#1}{\longrightarrow}}
\newcommand{\limx}[1]{\underset{#1}{\lim}}

\newcommand{\D}{\mathcal{D}}

\newcommand{\sotlimx}[1]{\text{SOT}\underset{#1}{\lim}}

\newcommand{\vaf}{\mathcal{V}(\A)_f}

\newcommand{\lwnorm}[1]{\bigl\Vert{#1} \bigr\Vert_{w,1}}

\newcommand{\lwnorma}[1]{\bigl\Vert{#1} \bigr\Vert_{L^{w,1}(\mathcal{A})}}

\newcommand{\lwnormb}[1]{\bigl\Vert{#1} \bigr\Vert_{L^{w,1}(\mathcal{B})}}

\newcommand{\intnorm}[1]{\bigl\Vert{#1} \bigr\Vert_{1\cap \infty}}
\newcommand{\inttnorm}[1]{\bigl\Vert{#1} \bigr\Vert_{L^{1}\cap L^{\infty}(\A)}}
\newcommand{\intxnorm}[2]{\bigl\Vert{#2} \bigr\Vert_{L^{1}\cap L^{\infty}(#1)}}

\newcommand{\pltnorm}[1]{\bigl\Vert{#1} \bigr\Vert_{L^{1}+ L^{\infty}(\A)}}
\newcommand{\plnnorm}[1]{\bigl\Vert{#1} \bigr\Vert_{L^{1}+ L^{\infty}(\nu)}}
\newcommand{\plxnorm}[2]{\bigl\Vert{#2} \bigr\Vert_{L^{1}+ L^{\infty}(#1)}}

\newcommand{\intnnorm}[1]{\bigl\Vert{#1} \bigr\Vert_{L^{1}\cap L^{\infty}(\nu)}}
\newcommand{\mtop}{\mathcal{T}_m}
\newcommand{\mtopc}{\overset{\mathcal{T}_m}{\rightarrow}}

\newcommand{\plc}{\overset{L^1+L^\infty}{\longrightarrow}}

\begin{document}

%
%
%
%
%
%
%
%
%

\title[Isometries on certain non-commutative spaces]{Extreme point methods in the study of isometries on certain non-commutative spaces}

\author[Pierre de Jager]{Pierre de Jager$^{*,1}$}

\address{Department of Mathematical Sciences\\ University of South Africa\\ P.O. Box, 392, Pretoria 0003, South Africa }\email{dejagp@unisa.ac.za}
\thanks{$^*$Corresponding author}

\thanks{$^1$ Department of Mathematical Sciences, University of South Africa, P.O. Box, 392, Pretoria 0003, South Africa, dejagp@unisa.ac.za}

\author[Jurie Conradie]{Jurie Conradie$^2$}
\address{Department of Mathematics \\ University of Cape Town \\ Cape Town, South Africa }
\email{ jurie.conradie@uct.ac.za}
\thanks{$^2$Department of Mathematics, University of Cape Town, jurie.conradie@uct.ac.za}


\date{29 July 2019}

\begin{abstract}
In this paper we characterize surjective isometries on certain classes of non-commutative spaces associated with semi-finite von Neumann algebras: the Lorentz spaces $L^{w,1}$, as well as the spaces $L^1+L^\infty$ and $L^1\cap L^\infty$. The technique used in all three cases relies on characterizations of the extreme points of the unit balls of these spaces. Of particular interest is that the representations of isometries obtained in this paper are global representations.
\end{abstract}

\maketitle

\section{Introduction}

If $T:X\to Y$ is an isometry of the normed space $X$ onto the normed space $Y$ and $x$ is an extreme point of the closed unit ball of $X$, $T(x)$ will be an extreme point of the closed unit ball of $Y$. This simple observation, when used in conjunction with characterizations of the extreme points of unit balls of certain classes of normed spaces, often plays an important role in determining the structure of isometries of these spaces. This is the  case in the characterization of isometries on the Lorentz spaces $L^{w,1}(0,1)$ in \cite{key-Car92} and on $L^1\cap L^\infty(0,\infty)$ and $L^1+L^\infty(0,\infty)$ in \cite{key-Grz92b}. Our aim in this paper is to derive characterizations of isometries on the non-commutative analogues of these spaces.

The spaces considered here are examples of non-commutative (quantum) symmetric spaces and consist of closed and densely defined operators on a Hilbert space $H$ affiliated with a von Neumann algebra $\A\subseteq B(H)$. We will confine ourselves to the setting where $\A$ is semi-finite and equipped with a distinguished faithful semi-finite normal trace and the operators are measurable with respect to this trace (precise definition will be given later). The special case where the trace of the identity is finite will be referred to as the \textit{finite setting}. If the von Neumann algebra is commutative and hence isometrically isomorphic to an $L^\infty$-space over some localizable measure space, one obtains the classical symmetric function spaces. 

The representations of isometries between various examples of symmetric spaces presented in the literature are typically in terms of a partial isometry or positive operator and a Jordan $*$-morphism between the underlying von Neumann algebras. Since a symmetric space is usually not contained in a von Neumann algebra, the representations are usually limited to the intersection of the von Neumann algebra and the symmetric space. By considering Jordan $*$-morphisms between spaces of trace-measurable operators, we are able to give representations valid on all of the symmetric space. We will refer to these as \textit{global representations}.

After a preliminary section on commutative and non-commutative function spaces, Jordan homomorphisms and extreme points, we consider surjective isometries between Lorentz spaces associated with semi-finite von Neumann algebras in Section \ref{S isometry}. A characterization of such isometries in the finite setting was obtained in \cite[Theorem 5.1]{key-Chilin89}.  In \cite[p.39]{key-Huang19} a structural description is given for into and onto isometries between more general Lorentz spaces but under the additional assumption that the isometries are positive. Our approach here is to use
 the structural description of a surjective isometry between a strongly symmetric space with order continuous norm and a symmetric space in our recent paper (see \cite[Theorem 5.3]{key-dJ19a}). Since this characterization requires  the isometry to be projection disjointness preserving and finiteness preserving, the main thrust of Section \ref{S isometry} is to show that a surjective isometry between the Lorentz spaces under consideration has these properties. Further structural analysis then allows for a refinement of the structural description provided by \cite[Theorem 5.3]{key-dJ19a}  and enables us to obtain a global representation. We also prove the converse of this result. 

In Section \ref{S4} we characterize surjective isometries between quantum $L^1+L^\infty$-spaces. The approach followed to obtain the desired structural description of such an isometry is to use a characterization of the extreme points of the unit ball of a $L^1+L^\infty$-space to show that the isometry restricted to the von Neumann algebra is an $L^\infty$-isometry. This leads to a local representation of the isometry in terms of a Jordan $*$-isomorphism of the von Neumann algebras using Kadison's description of isometries between $C^*$-algebras (see \cite[Theorem 7]{key-K51}). Further analysis enables one to extend the Jordan $*$-isomorphism to a Jordan $*$-isomorphism of the corresponding spaces of trace-measurable operators and obtain a global representation of the isometry. It is worth noting that since these spaces do not have order continuous norm, one does not have access to \cite[Theorem 5.3]{key-dJ19a}.

The final section is devoted to characterizing positive surjective isometries between non-commutative $L^1\cap L^\infty$-spaces. We use a description of the extreme points of the unit ball of $L^1\cap L^\infty$ to show that such an isometry is square-preserving and can therefore be extended to a Jordan $*$-homomorphism. Another possibility would have been to try to show that such isometries are finiteness preserving in order to apply \cite[Theorem 4.11]{key-dJ19a}. Since this would yield a less refined  structural description, we decided against this approach.

\section{Preliminaries}
\subsection{Function spaces}
We start by introducing the classical function spaces on which the non-commutative spaces considered in this paper are based. In what follows,  $(\Omega,\Sigma,\mu)$ will be an arbitrary $\sigma$-finite measure space. There are a number of different types of Lorentz spaces found in the literature (see, for example, \cite{key-BS88} and \cite{key-Car93}). In this paper we restrict our attention to the $L^{w,1}$-spaces defined in terms of a  decreasing weight function  $w:(0,\infty)\to (0,\infty)$  satisfying $\lim\limits_{t\to 0}w(t)=\infty, \lim\limits_{t\to\infty}w(t)=0, \int_0^\infty w(t)\,dt=\infty$ and $\int_0^1 w(t)\,dt=1$. Let $L_\infty^0(\Omega)$ denote the collection of all (equivalence classes of) measurable functions on $\Omega$ which are bounded except possibly on a set of finite measure. If $f\in L_\infty^0(\Omega)$, then the decreasing rearrangement $f^\ast$ exists (see \cite{key-BS88} for the definition and properties of decreasing rearrangements). We put $\norm{f}_{L^{w,1}(\Omega)}=\int_0^{\mu(\Omega)} f^\ast(t) w(t)\,dt$, and define the Lorentz function space $L^{w,1}(\Omega)$  as the set of all $f\in L_\infty^0(\Omega)$ for which $\norm{f}_{L^{w,1}(\Omega)}$ is finite, and  $\norm{\cdot}_{L^{w,1}(\Omega)}$ is  then a norm on $L^{w,1}(\Omega)$. If we define $\psi:[0,\infty)\to [0,\infty)$ by $\psi(t):=\int_0^t w(t)dt$, then these spaces are sometimes also denoted as $\Lambda_\psi(\mu)$, with norm $\norm{f}_{\Lambda_\psi(\mu)}:=\int_0^{\mu(\Omega)}f^*(t)d\psi(t)$. We note that since  $w$ is locally integrable, $\psi$ is continuous. Furthermore, $w(t)\neq 0$ for any $t\in [0,\tau(\id))$ and so $\psi$ is strictly increasing and therefore invertible, and  $\psi^{-1}$ is continuous.

The linear spaces $L^1\cap L^\infty (\Omega)=L^1(\Omega)\cap L^\infty (\Omega)$ and $L^1+ L^\infty (\Omega)=\{f\in L_\infty^0(\Omega):f=g+h, g\in L^1(\Omega), h\in L^\infty(\Omega)\}$ can be equipped with the norms 
$\intxnorm{\Omega}{f}=\text{max}\,\{\norm{f}_1,\norm{f}_\infty\}$ and $\plxnorm{\Omega}{f}=\text{inf}\,\{\norm{g}_1+\norm{h}_\infty:f=g+h, g\in L^1(\Omega), h\in L^\infty(\Omega)\}$, respectively. 

\subsection{The non-commutative setting}
 In order to define the non-commutative analogues of these function spaces we provide some background information regarding von Neumann algebras and trace-measurable operators. Suppose $H$ is a Hilbert space, $B(H)$ is the space of all bounded linear operators on $H$, and $\A\subseteq B(H)$ is a von Neumann algebra. We call $\A$ \emph{atomic} if it contains a set $\{p_\lambda\}_{\lambda \in \Lambda}$ of minimal projections such that $\sum p_\lambda =\id$ (where $\id$ denotes the identity operator on $H$), and this happens if and only if it is a product of Type 1 factors (see \cite[p.354]{key-Black06}). If $\mathcal{A}$ does not contain minimal projections, then it is called \emph{non-atomic}.  We will use the phrase ``($\A$, $\tau$) is a semi-finite von Neumann algebra'' to describe the situation where $\A$ is a semi-finite von Neumann algebra equipped with a distinguished faithful normal semi-finite trace $\tau$ and unless stated otherwise this will be the setting in what follows. Let $\pa$ denote the lattice of projections in $\A$ and $\paf$ the sublattice consisting of projections with finite trace. Two projections $p,q\in\pa$ are \emph{equivalent}, written $p\sim q$, if there is a $v\in\A$ such that $v^\ast  v=p$ and $vv^\ast=q$. If $p\sim q$, then $\tau(p)=\tau(q)$. A projection $p\in\pa$ is called a \emph{finite projection} if there is no proper subprojection of $p$ equivalent to $p$. Any projection with finite trace is finite, but the converse need not be true. We will use the notation $\rax{\A}$, $\rax{SOT}$ and $\rax{WOT}$ to denote convergence in $\A$ with respect to the  operator norm topology, the strong operator topology (SOT) and the weak operator topology (WOT), respectively. Further details regarding von Neumann algebras may be found in \cite{key-K1}.

Let $x$ be a closed densely defined operator on $H$ with domain $\D(x)$. We will use $n(x)$, $r(x)$ and $s(x)$ to denote the projection onto the kernel of $x$, the projection onto the closure of the range of $x$ and the \textit{support projection} $\id - n(x)$, respectively. We say that $x$ is \textit{affiliated} with $\A$ if $u^*xu=x$ for all unitary operators $u$ in the commutant of $\A$. If, in addition, there is a sequence $(p_n)_{n=1}^\infty$ in $\pa$ such that $p_n \uparrow \id$, $p_n(H)\subseteq \D(x)$ and $\id - p_n\in\paf$ for every $n$, then $x$ is called \textit{$\tau$-measurable}.  The set of all $\tau$-measurable operators will be denoted by $\tma$. Equipped with the measure topology $\mtop$ (see \cite[p.210]{key-Dodds14}), it becomes a complete metrisable topological $\ast$-algebra when sums and products are defined as the closures of the algebraic sum and algebraic product, respectively. We can define a partial order on the self-adjoint elements of $\tma$ in the following way. Let $\langle\cdot,\cdot\rangle$ denote the inner product on $H$. We will say that $x\in \tma$ is positive if  $\langle x\eta,\eta\rangle\geq 0$ for all $\xi$ in the domain of $x$, and we use the notation $\tma^+=\{x\in\tma: x\geq 0\}$. A partial order can then be defined on $\tma^{sa}=\{x\in \tma:x=x^*\}$ using the cone $\tma^+$.  We will write $x_\lambda \uparrow x$ if $(x_\lambda)_{\lambda \in \Lambda}$ is an increasing net in $\tma$ with supremum $x$. Note that for any $x\in \tma^+$ there exists an increasing net $\net{x}$ of positive elements in $\mathcal{F}(\tau)=\{x\in \A:\tau(s(x))<\infty\}$ such that $x_\lambda \uparrow x$ (see \cite[Proposition 1 vii)]{key-Dodds14}). For more information about $\tau$-measurable operators the  reader is referred to \cite{key-Dodds14} and \cite{key-Terp1}. 

To each $x\in \tma$ we can associate a positive decreasing function $\svf{x}$ on $(0,\infty)$ called the \emph{singular value function} of $x$. For $t>0$ this function is defined as $\svft{x}{t}=\inf\{s\geq 0: \dfs{|x|}{s}\leq t\}$, where $\dfs{|x|}{s}:=\tau\left(e^{|x|}(s,\infty)\right)$ ($s\geq 0$) is the \emph{distribution function} of $|x|$. The singular value function can be used to extend the trace $\tau$ to an additive, positively homogeneous, unitarily invariant normal map on $\tma^+$ by setting $\tau(x):=\int_0^{\infty}\svft{x}{t}dt$ for $x\in \tma^+$. Suppose $E\subseteq S(\mathcal{A},\tau)$ is a linear subspace equipped with a norm $\norm{\cdot}_E$.  If $E$ is a Banach space and $x\in E$ with $\norm{x}_E\leq \norm{y}_E$, whenever $y\in E$, $x\in S(\mathcal{A},\tau)$ and $\svf{x}\leq \svf{y}$, then $E$ is called a \emph{symmetric space}. If this is the case and $x\in E$, then $x$, $x^*$ and $|x|$ all have the same norm and $uxv\in E$ with $\norm{uxv}_E\leq \norm{u}_{\mathcal{A}}\norm{v}_{\mathcal{A}}\norm{x}_E$ whenever $u,v\in \mathcal{A}$. To describe how one can obtain non-commutative analogues of classical function spaces we let $L^\infty(0,\infty)$ denote the collection of all essentially bounded Lebesgue measurable functions on $(0,\infty)$. If $\nu$ is given by integration with respect to Lebesgue measure, then $\nu$ yields a faithful normal semi-finite trace on the commutative von Neumann algebra $L^\infty(0,\infty)$. In this case $S(L^\infty(0,\infty),\nu)=L_\infty^0(0,\infty)$ is the space of all Lebesgue measurable functions on $(0,\infty)$ that are bounded except possibly on a set of finite measure and the singular value function corresponds to the decreasing rearrangement of a measurable function. Suppose $(\mathcal{A},\tau)$ is a semi-finite von Neumann algebra and $E(0,\infty)\subseteq L_\infty^0(0,\infty)$ is a symmetric space. Letting $E(\A):=\{x\in S(\mathcal{A},\tau):\svf{x}\in E(0,\infty)\}$ we obtain a linear subspace of $\tma$ and one can define a norm on $E(\A)$ using $\norm{x}_{E(\A)}=\norm{\svf{x}}_{E(0,\infty)}$, for $x\in E(\A)$. It can be shown that $E(\A)$,  equipped with this norm, is a symmetric space (see \cite{key-Kalton08}).  In particular, we can define the non-commutative spaces $L^{w,1}(\A)$, $L^1\cap L^\infty(\A)$ and $L^1+L^\infty(\A)$ in this way. It can be shown (see \cite[p.725]{key-Dodds93}) that $L^1\cap L^\infty(\A)=L^1(\A)\cap \A$, $L^1+L^\infty(\A)=L^1(\A)+\A$, $\inttnorm{x}=\max\{\norm{x}_1,\norm{x}_\A\}$ and $\pltnorm{x}=\text{inf}\,\{\norm{y}_1+\norm{z}_\infty:x=y+z, y\in L^1(\A), z\in \A\}$. To simplify notation we will typically denote the norms of $L^{w,1}(\A)$, $L^1\cap L^\infty(\A)$ and $L^1+L^\infty(\A)$ by $\norm{\cdot}_{w,1}$, $\norm{\cdot}_{1\cap \infty}$ and $\norm{\cdot}_{1+\infty}$, respectively. These three spaces are, in fact, \emph{fully symmetric spaces}, since if $E$ is one of these spaces, then it satisfies the additional property that $x\in E$ and $\norm{x}_E \leq \norm{y}_E$ whenever $x\in S(\mathcal{A},\tau)$, $y \in E$ and $x \prec\prec y$ (we write $x \prec\prec y$ if $\int_0^t \svft{x}{s}ds\leq \int_0^t \svft{y}{s}ds$ for all $t> 0$).  Further information about symmetric spaces may be found in \cite{key-Dodds14} and \cite{key-Pag}.

\begin{rem} \label{RT2 29/01/15}
The norm $\norm{\cdot}_E$ on a symmetric space $E$ is called order continuous if $\norm{x_\lambda}\downarrow 0$ whenever $x_\lambda \downarrow 0$ in $E$. Since $L^{w,1}(0,\infty)$ has order continuous norm (see \cite[Corollary 1]{key-For12b}), it follows from \cite[Theorem 54]{key-Dodds14} that $L^{w,1}(\A)$ has order continuous norm. Therefore, $\mathcal{F}(\tau)$ is norm dense in $L^{w,1}(\A)$.  It can be shown, using the spectral theorem, that for every $x\in \mathcal{F}(\tau)$, there is a sequence $(x_n)_{n=1}^\infty$ in $\mathcal{G}(\A)_f$ (the set of all finite linear combinations of mutually orthogonal projections in $\paf$) such that $x_n\rax{\A} x$ and also $x_n\rax{E} x$ (see \cite[Proposition 2.1]{key-dJ19a}). It follows that $\mathcal{G}(\A)_f$ is dense in $L^{w,1}(\A)$ and $\mathcal{G}(\A)_f^+$ is dense in $L^{w,1}(\A)^+$, using the order continuity of the norm.
\end{rem}

Suppose $(\A\subseteq B(H),\tau)$ and $(\B\subseteq B(K),\nu)$ are semi-finite von Neumann algebras, $\mathcal{H}\subseteq \tma$ is a linear subspace and $T:\mathcal{H}\ra \nmb$ is a linear map. If $\paf \subseteq \mathcal{H}$, then $T$  will be called \emph{finiteness preserving} if $\nu(s(T(p)))<\infty$ for every $p\in \paf$ and \emph{projection disjointness preserving} if $p,q\in \paf$ with $pq=0$ implies that $T(p)^*T(q)=0=T(p)T(q)^*$. The map $T$ will be called \emph{normal} (on $\mathcal{H}$) if $T(x_\lambda)\uparrow T(x)$ whenever $\{x_\lambda\}_{\lambda \in \Lambda}\subseteq\mathcal{H}^{sa}$ is an increasing net with supremum $x\in \mathcal{H}^{sa}$. 

\subsection{Jordan homomorphisms}

Suppose $\Phi:\A\to \B$ is a linear map. If $\Phi(yx + xy) = \Phi(y)\Phi(x) + \Phi(x)\Phi(y)$ for all $x, y \in \A$, then $\Phi$ is called a \textit{Jordan homomorphism}. We call such a map a \textit{Jordan $\ast$-homomorphism}, if, in addition, it also preserves adjoints. Analogously, we can define Jordan homomorphisms/$*$-homomorphisms between $\ast$-algebras of trace-measurable operators and these will play a significant role in global representations of isometries. We will need the following extension result.

\begin{lem}\label{L2 17/09/19}
Suppose $(\A,\tau)$ and $(\B,\nu)$ are semi-finite von Neumann algebras and $\Phi:\A \ra \B$ is a Jordan $*$-isomorphism. If $\nu(\Phi(p))=\tau(p)$ for every $p\in \paf$, then $\Phi$ has a unique extension to a Jordan $*$-isomorphism $\tilde{\Phi}$ from $\tma$ onto $\nmb$ and $\tilde{\Phi}$ is trace-preserving on $\tma$. 
\begin{proof}
It is clear that $\nu(\Phi(x))=\tau(x)$, whenever  $x\in \mathcal{G}(\A)_f$. If $x\in \A^+$, then there exists an increasing net $\net{x}\subseteq \mathcal{G}(\A)_f^+$ such that $x_\lambda \uparrow x$. Since $\Phi$ is a Jordan $*$-isomorphism it is an order-isomorphism and hence normal. Using the normality of $\tau$, $\nu$ and $\Phi$, we obtain $\tau(x_\lambda)\uparrow \tau(x)$ and $\nu(\Phi(x_\lambda))\uparrow \nu(\Phi(x))$ and hence $\nu(\Phi(x))=\tau(x)$.  Since any element in $\A$ can be written as a linear combination of positive elements, we have that $\Phi$ is trace-preserving. Therefore $\nu\circ\Phi$ is $\epsilon - \delta$ absolutely continuous with respect to $\tau$ on $P(\A)$ (i.e. for any $\epsilon>0$ there exists a $\delta>0$ such that $\nu(\Phi(p))<\epsilon$ whenever $p\in P(\A)$ is such that $\tau(p)<\delta$). It follows by \cite[Proposition 4.1]{key-Lab13} that  $\Phi$ extends uniquely to a continuous (with respect to the measure topologies) Jordan $*$-homomorphism $\tilde{\Phi}:\tma \ra \nmb$. Since it is easily checked that $\Phi^{-1}$ is also trace-preserving, one can similarly extend $\Phi^{-1}$ to a continuous (with respect to the measure topologies) Jordan $*$-homomorphism $\widetilde{\Phi^{-1}}:\nmb \ra \tma$. We show that $\tilde{\Phi}$ is injective. Suppose $x,y\in \tma$ with $\tilde{\Phi}(x)=\tilde{\Phi}(y)$. Using the density (with respect to the measure topology) of $\A$ in $\tma$ (see \cite[p. 210]{key-Dodds14}), there exist sequences $\seq{x},\seq{y}$ in $\A$ such that $x_n\mtopc x$ and $y_n \mtopc y$. Therefore  $\Phi(x_n)=\tilde{\Phi}(x_n)\mtopc \tilde{\Phi}(x)$ and $\Phi(y_n)=\tilde{\Phi}(y_n)\mtopc \tilde{\Phi}(y)$. Since $\widetilde{\Phi^{-1}}$ and $\Phi^{-1}$ agree on $\B$, we have that $x_n=\widetilde{\Phi^{-1}}(\Phi(x_n)) \mtopc \widetilde{\Phi^{-1}}(\tilde{\Phi}(x))$ and $y_n=\widetilde{\Phi^{-1}}(\Phi(y_n)) \mtopc \widetilde{\Phi^{-1}}(\tilde{\Phi}(y))=\widetilde{\Phi^{-1}}(\tilde{\Phi}(x))$. Therefore $x=y$. Similar arguments can be used to show that $\tilde{\Phi}$ is surjective and $\tilde{\Phi}^{-1}=\widetilde{\Phi^{-1}}$. Since $\tilde{\Phi}$ and $\tilde{\Phi}^{-1}$ are positive (see \cite[Lemma 3.9]{key-Weigt09}), it is easily verified that $\tilde{\Phi}$ is normal. Furthermore, for any $x\in \tma^+$ there exists an increasing net $\net{x}\subseteq \mathcal{F}(\tau)^+$ such that $x_\lambda \uparrow x$ and therefore a similar argument can be used to show that $\tilde{\Phi}$ is trace-preserving on $\tma$.
\end{proof}
\end{lem}

\subsection{Extreme points}\label{S extreme}

For a normed space $E$, let $B_E:=\{x\in E:\norm{x}\leq 1\}$ denote the unit ball in $E$ and $S_E:=\{x\in E:\norm{x}= 1\}$ the unit sphere in $E$. We wish to obtain characterizations of the extreme points of the unit balls of $L^{w,1}(\A)$, $L^1+ L^{\infty}(\A)$ and $L^1\cap L^\infty(\A)$. When dealing with isometries of $L^1\cap L^{\infty}(\A)$ and $L^1+L^\infty(\A)$ our techniques will restrict us to the non-atomic setting. In this setting the desired characterizations of the extreme points can be deduced from the corresponding results in the commutative setting (see \cite[Proposition 2.2]{key-Car92}, \cite[Proposition 2.3]{key-Schaefer92} and \cite[Corollary 1]{key-Grz92}, respectively) using the following result.

\begin{thm}$(\text{cf. \cite[Theorem 1.1, Corollary 5.18 and the Remark on p.223]{key-Chilin92}})$ \label{T2.2 Czer12} 
Suppose $(\A,\tau)$ is a non-atomic semi-finite von Neumann algebra and $E(0,\tau(\id))$ is a fully symmetric space. An operator $x\in S_{E(\A)}$ is an extreme point of $B_{E(\A)}$ if and only if $\svf{x}\restriction_{(0,\tau(\id))}$ is an extreme point of $B_{E(0,\tau(\id))}$ and one of the following conditions holds:
\begin{enumerate}
\item $\limx{t\ra \infty} \svft{x}{t}=0$
\item $n(x)\A n(x^*)=\{0\}$ and $|x|\geq \limx{t\ra \infty} \svft{x}{t}s(x)$.
\end{enumerate}
\end{thm}

\begin{rem}\label{rT2.2 Czer12} 
In \cite[Theorem 2.2, p.37]{key-Czer12} it is claimed that the restriction to non-atomic von Neumann algebras can be removed by embedding an arbitrary semi-finite von Neumann algebra into a non-atomic one. The image of the original algebra in the non-atomic algebra need however not be a non-atomic subalgebra. The following example shows that Theorem \ref{T2.2 Czer12} may fail if the original algebra is atomic. Let $H=\ell^2(\mathbb{N})=\ell^2$ and $\A=\{M_f:f\in \ell^\infty\}$. Then $\A$ is isometrically $*$-isomorphic to $\ell^\infty$ and $L^1(\A)$ is isometrically $*$-isomorphic to $\ell^1$. Let $(e_{i,n})_{n=1}^\infty$ be the sequence in $\ell^1$ defined by $e_{i,n}=\delta_{i}(n)$. Then for each $i\in \mathbb{N}$, $(e_{i,n})_{n=1}^\infty$ is an extreme point of $B_{\ell^1}$ (or one could use \cite[Theorem 3]{key-For10} to conclude that the set of extreme points of $B_{\ell^1}$ is non-empty). It is however the case that $B_{L^1(0,\infty)}$ does not have extreme points (see \cite[Remark 1]{key-For10}) and in particular $\svf{(e_{i,n})_{n=1}^\infty}=\chi_{[0,1)}$ is not an extreme point of $B_{L^1(0,\infty)}$. 
\end{rem}

In order to characterize the extreme points of the unit ball of $L^{w,1}(\A)$ and to facilitate the application of Theorem \ref{T2.2 Czer12} we will repeatedly use the following two lemmas.

\begin{lem} \label{LP3 24/10/14}
If $x\in \tma$ is such that $\tau(r(x))<\infty$ (or equivalently $\tau(s(x))<\infty$), then there exists a unitary $u\in \A$ such that $x=u|x|$.
\begin{proof}
We start by noting that since $r(x)$ and $s(x)$ are equivalent projections, $\tau(r(x))=\tau(s(x))$,  and hence $\tau(r(x))<\infty$ is equivalent to $\tau(s(x))<\infty$. Let $x=v|x|$ be the polar decomposition of $x$. Let $p=v^*v$ and $q=vv^*$. Then $p\sim q$ and hence $\tau(p)=\tau(q)=\tau(r(x))<\infty$. Since any projection with finite trace is finite, $p$ and $q$ are finite projections. By \cite[Exercise 6.9.7(6)]{key-K2}, $p^{\perp}\sim\, q^{\perp}$. We can therefore find a partial isometry $w\in \A$ such that $w^*w=p^{\perp}$ and $ww^*=q^{\perp}$. Note that, using \cite[Exercise 2.8.45]{key-K1}, we obtain $r(w)=q^{\perp}=r(x)^{\perp}=r(v)^{\perp}=n(v^*)$ and hence $v^*w=0$. We can similarly show that $vw^*=0$ and hence $w^*v=0=wv^*$. Let $u=v+w$. Then $u^*u=v^*v+v^*w+w^*v+w^*w=p+0+0+p^{\perp}=\id$. Similarly $uu^*=\id$, and hence $u$ is unitary. Furthermore, $w|x|=ws(w)r(|x|)|x|=wp^\perp p|x| =0$, by \cite[Exercise 2.8.45]{key-K1} and using the facts that $r(|x|)=r(x^*)=s(x)=p$ and $s(w)=p^{\perp}$. It follows that $u|x|=(v+w)|x|=v|x|+0=x$.
\end{proof}
\end{lem}

\begin{lem} \label{LP1 20/05/15}
Suppose $x\in \tma$ and $\beta>0$. Then $\svf{x}=\beta\chi_{[0,\alpha)}$ if and only if $|x|=\beta p$ for some $p\in \pa$ with $\tau(p)=\alpha$.
\begin{proof}
If $|x|=\beta p$ for some $p\in \pa$, then a direct calculation using the definition of the singular value function shows that $\svf{x}=\beta \chi_{[0,\tau(p))}$. Conversely, suppose $\svf{x}= \beta\chi_{[0,\alpha)}$. It is easily checked that $\tau(e^{|x|}(s,\infty))=\alpha$ for all $s<\beta$. Using the faithfulness of $\tau$ this implies that $e^{|x|}(0,s]=0$ for all $s<\beta$. Furthermore, $e^{|x|}(s,\infty)=0$ for all $s>\beta$, since $\norm{x}_{\A}=\beta$. It follows that $|x|$ has two eigenvalues, namely $0$ and $\beta$, and therefore $|x|=\beta p$ for some $p\in \pa$ and $\tau(p)=\alpha$.
\end{proof}
\end{lem}

Using Lemma \ref{LP3 24/10/14} we obtain the following semi-finite extension of \cite[Lemma 2.3]{key-Chilin89}.

\begin{lem} \label{L2.3 Chilin1} Suppose $E(\tau)$ is a symmetric space and $x\in B_{E(\tau)}$ is such that $\tau(r(x))<\infty$ (or equivalently $\tau(s(x))<\infty$). Then $x$ is an extreme point of $B_{E(\tau)}$ if and only if $|x|$ is one.
\end{lem}

\section{Surjective isometries between Lorentz spaces}\label{S isometry}

Throughout this section we will assume that $(\A,\tau)$ and $(\B, \nu)$ are semi-finite von Neumann algebras and that $w$ is a strictly decreasing weight function. The aim of this section is to characterize surjective isometries from $L^{w,1}(\A)$ onto $L^{w,1}(\B)$. That $w$ is strictly decreasing will be used in the characterization of the extreme points of the unit ball of $L^{w,1}(\A)$. Furthermore, it follows from $\int_0^\infty w(t)\,dt=\infty$ that $w(t)>0$ for all $t\geq 0$ and hence enables one to show  that $L^{w,1}(\A)$ has strictly monotone norm (see \cite[p.532]{key-Chilin89}, for example). Suppose $U:L^{w,1}(\A)\ra L^{w,1}(\B)$ is a surjective isometry. In the setting where $\tau(\id)<\infty$ and $\nu(\id)<\infty$ (see \cite[Theorem 5.1]{key-Chilin89}), which we will refer to as the \emph{finite setting}, the representation of $U$ is obtained by using the characterization of the extreme points of the unit ball of a Lorentz space to show that $U(\id)=\frac{1}{\psi(\nu(|a|))}a$ for some partial isometry $a\in \B$. One can then find a unitary operator $u\in \B$ such that $a=u|a|$. The most substantial part of the proof involves showing that the surjective isometry $T(x):=u^*U(x)$ is positive and that $|a|=\id$. A structural description of positive surjective isometries between a symmetric space and fully symmetric space (\cite[Theorem 3.1]{key-Chilin89}) is then employed to obtain the desired representation. Recently it has been shown that \cite[Theorem 3.1]{key-Chilin89} can be extended to the semi-finite setting and that this extension can then be used to obtain a structural description of projection disjointness and finiteness preserving surjective isometries between symmetric spaces associated with semi-finite von Neumann algebras (see \cite[Theorem 5.3]{key-dJ19a}). We will show that applying the techniques employed in the finite setting to $U(p)$, for each $p\in \paf$, (instead of $U(\id)$, which need not be defined in the semi-finite setting)  will enable us to show that $U$ is disjointness-preserving in the semi-finite setting. Furthermore, the characterization of the extreme points of the unit ball of a Lorentz space will be used to show that $U$ is finiteness-preserving. This will enable us to use \cite[Theorem 5.3]{key-dJ19a} to obtain a preliminary structural description of the isometry $U$. Further analysis of this structural description in the context of Lorentz spaces will yield the desired representation. For the sake of brevity we have omitted those portions of the proof that involve only minor modifications to arguments in the finite setting or which use very natural non-commutative analogues of arguments employed in the commutative setting. The interested reader is referred to the first author's doctoral thesis (\cite[$\S$7.2]{key-dJ17}) for a detailed verification of these components of the proof.

We start by extending the characterization of the extreme points of the unit ball of $L^{w,1}(\A)$ to the semi-finite setting. We need one preliminary result. By replacing decreasing rearrangements with singular value functions in the proof of \cite[Lemma 3]{key-For12b} and restricting to the situation where $w$ is strictly decreasing and $\int_0^\infty w(t)dt=\infty$ in \cite[Lemma 3]{key-For12b}, we obtain the following non-commutative analogue of this result.  

\begin{lem} \label{L2.1 Car} Suppose $(\A ,\tau)$ is a semi-finite von Neumann algebra, $w:(0,\infty)\ra (0,\infty)$ is a strictly decreasing weight function and $x=\alpha p\in S_{L^{w,1}(\A)}$, for some $p\in \paf$ and $\alpha=1/\int_0^{\tau(p)}w(t)dt$.  If there exist $y,z\in S_{L^{w,1}(\A)}$ with $x=\tfrac{1}{2}(y+z)$, then $\svf{x}=\tfrac{1}{2}(\svf{y}+\svf{z})$.
\end{lem}

Since we are working with arbitrary semi-finite von Neumann algebras (not just non-atomic ones), we cannot deduce the structure of extreme points from the commutative setting (see Remark \ref{rT2.2 Czer12}). The sufficiency part of the following result can, however, be obtained using a non-commutative analogue of the corresponding argument in the proof of \cite[Proposition 2.2]{key-Car92}; and the necessity by adapting the corresponding part of the proof of \cite[Theorem 4.1]{key-Chilin89} to the semi-finite setting.

\begin{prop} \label{AnP2.2 Car} 
Suppose $(\A,\tau)$ is a semi-finite von Neumann algebra, $w$ is a strictly decreasing weight function and $E=L^{w,1}(\A)$. Then $x$ is an extreme point of $B_E$ if and only if $x=\frac{v}{\psi(\tau(|v|))}$ for some $v\in \vaf$, where $\vaf$ denotes the set of partial isometries $w$ in $(\A,\tau)$ with $\tau(|w|)<\infty$.
\end{prop}

Suppose $U:L^{w,1}(\A) \ra L^{w,1}(\B)$ is a surjective isometry and $0\neq p\in \paf$. Then $x=\frac{1}{\psi(\tau(p))}p$ is an extreme point of the unit ball of $L^{w,1}(\A)$, by Proposition \ref{AnP2.2 Car}. Since $U$ is a surjective isometry, $U(x)$ is also an extreme point and therefore $U(p)$ can be written in the form $\alpha_p v_{p}$, where $v_{p} \in \mathcal{V}(\B)^f$ and $\alpha_p = \frac{\psi(\tau(p))}{\psi(\nu(|v_{p}|))}$. By Lemma \ref{LP3 24/10/14}, there exists a unitary operator $u_{p}\in \B$ such that $v_{p}=u_{p}|v_{p}|$. Define $T_p:L^{w,1}(\A)\ra L^{w,1}(\B)$ by $T_p(x):=u_{p}^*U(x)$ for $x\in L^{w,1}(\A)$. It is easily checked that $T_p$ is a surjective isometry. 

We wish to show that $T_p(pL^{w,1}(\A)p)\subseteq |v_{p}|L^{w,1}(\B)|v_{p}|$ and $T_p(x)\geq 0$ for every $x\in pL^{w,1}(\A)^+p$. Let $0\neq q\in \paf$ with $q\leq p$.  Since $T_p$ is a surjective isometry, we have that $T_p(q)=\beta v$ for some $v\in \mathcal{V}(\B)^f$, where $\beta = \frac{\psi(\tau(q))}{\psi(\nu(|v|))}$. Let $y=T_p(p+q)$. Then $y=u_{p}^*U(p)+T_p(q) = \alpha_p |v_{p}|+\beta v$. It follows that $s(y)\leq s(|v_{p}|)\vee s(v)=|v_{p}|\vee |v|$. Furthermore, $v_{p}, v\in \mathcal{V}(\B)^f$ and so $\nu(|v_{p}|),\nu(|v|)<\infty$. It follows that $|v_{p}|\vee |v|$ (and hence $s(y)$) has finite trace. By Lemma \ref{LP3 24/10/14} this implies that there exists a unitary operator $w\in \B$ such that $wy=|y|$. By making suitable adjustments to the arguments employed in the proofs of \cite[Lemmas 5.2, 5.3, 5.4 and 5.6]{key-Chilin89} we obtain $w|v_{p}|=|v_{p}|$, $wv=|v|$ and $|v|\leq |v_{p}|$.

\begin{lem} \label{L5.7 C1}
$T_p(pL^{w,1}(\A)p)\subseteq |v_{p}|L^{w,1}(\B)|v_{p}|$.
\begin{proof}
Since $|v|\leq |v_{p}|$, we have $|v_{p}||v||v_{p}|=|v|$ and so $|v|\in |v_{p}|L^{w,1}(\B)|v_{p}|$. Therefore $v\in |v_{p}|L^{w,1}(\B)|v_{p}|$ and so $T_p(q)=\beta v\in |v_{p}|L^{w,1}(\B)|v_{p}|$. Let $\mathcal{D}_{p}:=\{e\in \paf :e\leq p\}$ and $\mathcal{G}_p^f:=\text{span}\,(\mathcal{D}_{p})$. Since $0<q\leq p$ was arbitrary, we have that $ T_p(\mathcal{G}_p^f)\subseteq |v_{p}|L^{w,1}(\B)|v_{p}|$. $T_p$ is continuous and it is easily checked that $|v_{p}|L^{w,1}(\B)|v_{p}|$ is closed in $L^{w,1}(\B)$; therefore, $T_p(\overline{\mathcal{G}_p^f})\subseteq |v_{p}|L^{w,1}(\B)|v_{p}|$. This completes the proof of the lemma, since one can use the density results given in Remark \ref{RT2 29/01/15} to show that  $\overline{\mathcal{G}_p^f}=pL^{w,1}(\A)p$. 
\end{proof}
\end{lem} 

\begin{lem} \label{L5.8 C1}
If $x\in pL^{w,1}(\A)^+p$, then $T_p(x)\geq 0$.
\begin{proof}
Recall that $q\in \paf$ with $0\neq q \leq p$ and that $T_p(q)=\beta v$.  Since $v\in |v_{p}|L^{w,1}(\B)|v_{p}|$, this implies that $\beta v =\beta |v_{p}|v$. Furthermore, $w|v_{p}|=|v_{p}|$, $wv=|v|$ and it can be shown (as in the proof of \cite[Lemma 5.4]{key-Chilin89}) that $w$ and $|v_{p}|$ commute. We therefore obtain  \[T_p(q)=\beta |v_{p}|v=\beta w|v_{p}|v= \beta |v_{p}|wv =\beta |v_{p}||v|.\] Since $|v|\leq |v_{p}|$, it follows that $|v_{p}||v|=|v|$ and hence $T_p(q)=\beta |v|\geq 0$.  It is then easily checked that $T_p(\mathcal{G}_p^+)\subseteq  L^{w,1}(\B)^+$, since $0<q\leq p$ was arbitrary. Furthermore, $\overline{\mathcal{G}_p^+}= pL^{w,1}(\A)^+p$ (see Remark \ref{RT2 29/01/15}), $T_p$ is continuous and $ L^{w,1}(\B)^+$ is closed (see \cite[Corollary 12]{key-Dodds14}). It follows that $T_p(pL^{w,1}(\A)^+p)=T_p(\overline{\mathcal{G}_p^+})\subseteq L^{w,1}(\B)^+$.
\end{proof}
\end{lem}

The major part of the proof up to this point has been a semi-finite adaptation of the techniques employed in the proof of \cite[Theorem 5.1]{key-Chilin89}. We show that this groundwork in fact enables us to prove that any surjective isometry between Lorentz spaces is projection disjointness preserving. We will need the following easily verified claim.

\begin{lem} \label{LP10 26/09/14} 
If $0\neq p,q,e\in \pa$ and $\alpha, \beta, \gamma >0$, then $\alpha p+\beta q =\gamma e$ if and only if one of the following conditions holds
\begin{enumerate}
\item$p=q=e$ and $\alpha+\beta =\gamma$ or
\item $pq=0$, $p+q=e$ and  $\alpha=\beta=\gamma$.
\end{enumerate}
\end{lem}

\begin{prop} \label{TL5.9 C1}\label{TC5.0 C1} 
Suppose $(\A,\tau)$ and $(\B, \nu)$ are semi-finite von Neumann algebras and $w$ is a strictly decreasing weight function. If $U:L^{w,1}(\A)\ra L^{w,1}(\B)$ is a surjective isometry, then $U(p)^*U(q)=0=U(p)U(q)^*$ whenever $p,q\in \paf$ are such that $pq=0$. Furthermore,  $\alpha_p=\alpha_q=\alpha_{p+q}$, $v_{p}+v_{q}=v_{p+q}$ and $|v_{p}|+|v_{q}|=|v_{p+q}|$, where $v_{p},v_{q}, v_{p+q}\in \mathcal{V}(\B)^f$ and $\alpha_p,\alpha_q,\alpha_{p+q}\in \mathbb{R}$ are such that $U(p)=\alpha_p v_{p}$, $U(q)=\alpha_q v_{q}$ and $U(p+q)=\alpha_{p+q}v_{p+q}$. 
\begin{proof}
Since $p+q=p \vee q \in \paf$, there exits a partial isometry $v_{p+q}\in \mathcal{V}(\B)^f$ such that $U(p+q)=\alpha_{p+q}v_{p+q}$. By Lemma \ref{LP3 24/10/14}, there exists a unitary operator $u_{p+q}\in \B$ such that $v_{p+q}=u_{p+q}|v_{p+q}|$. Let $u$ denote $u_{p+q}$ and note that $T_{p+q}(p)=u^*U(p)= \alpha_p u^*v_{p}$. Furthermore, $(u^*v_{p} )^*(u^*v_{p} )=v_{p}^*v_{p}=|v_{p}|$. It follows by \cite[Proposition 6.1.1]{key-K2} that $\frac{1}{\alpha_p}T_{p+q}(p)$ is a partial isometry. Since $\frac{1}{\alpha_p}T_{p+q}(p)\geq 0$, by Lemma \ref{L5.8 C1}, $\frac{1}{\alpha_p}T_{p+q}(p)$ is a projection. We can similarly show that $\frac{1}{\alpha_q}T_{p+q}(q)$ and $\frac{1}{\alpha_{p+q}}T_{p+q}(p+q)$ are also projections. Therefore,  $u^*v_{p}$, $u^*v_{q}$ and $u^*v_{p+q}$ are all projections. We also have 
\begin{eqnarray*}
&\alpha_p u^*v_{p}+\alpha_q u^*v_{q}=T_{p+q}(p)+T_{p+q}(q) = T_{p+q}(p+q)=\alpha_{p+q} u^*v_{p+q}.&
\end{eqnarray*}
By Lemma \ref{LP10 26/09/14} this implies that 
\begin{eqnarray}
&u^*v_{p}=u^*v_{q}=u^*v_{p+q} \,\, \text{and}\,\, \alpha_p+\alpha_q=\alpha_{p+q} \,\, \text{or}&  \label{e5.9.0 C1} \\
&(u^*v_{p})(u^*v_{q})=0, u^*v_{p}+u^*v_{q}=u^*v_{p+q}\,\, \text{and}\,\,  \alpha_p=\alpha_q=\alpha_{p+q}.& \label{e5.9.1 C1}
\end{eqnarray}
If we assume that (\ref{e5.9.0 C1}) holds, then $T_{p+q}\left(\frac{p}{\alpha_p}\right)=u^*v_{p}=u^*v_{q}=T_{p+q}\left(\frac{q}{\alpha_p}\right)$. Since $T_{p+q}$ is an isometry and hence injective, it follows that $\frac{p}{\alpha_p}=\frac{q}{\alpha_q}$ and hence $p=q$. This is a contradiction, since $pq=0$ and $p,q\neq 0$. Therefore (\ref{e5.9.1 C1}) holds. Note that $v_{q}^*u=(u^*v_{q})^*=u^*v_{q}$, since $u^*v_{q}$ is a projection. Similarly, $v_{p}^*u=u^*v_{p}$.  Therefore  $(u^*v_{p})(v_{q}^*u) =(u^*v_{p})(u^*v_{q})=0$, using (\ref{e5.9.1 C1}). It follows that $v_{p} v_{q}^*=0$. Furthermore, $(u^*v_{p})(u^*v_{q})=0$ implies that $v_{p}^*v_{q}=  (v_{p}^*u)(u^*v_{q})=0$, since $uu^*=\id$ and $v_{p}^*u=u^*v_{p}$. Therefore $U(p)^*U(q)=0=U(p)U(q)^*$.

 Note also that $u(u^*v_{p}+u^*v_{q})=u(u^*v_{p+q})$, using (\ref{e5.9.1 C1}). Since $u$ is unitary, it follows that $v_{p}+v_{q}=v_{p+q}$. Furthermore, $v_{p}^*v_{q}=0=v_{q}^*v_{p}$ and so $|v_{p+q}|=v_{p+q}^*v_{p+q}=(v_{p}+v_{q})^*(v_{p}+v_{q})=v_{p}^*v_{p}+v_{q}^*v_{q}=|v_{p}|+|v_{q}|$.
\end{proof}
\end{prop}

We are now in a position to prove the main result of this section.  

\begin{thm} \label{T5 C1} Suppose $(\A,\tau)$ and $(\B, \nu)$ are semi-finite von Neumann algebras and $w$ is a strictly decreasing weight function.
\begin{itemize}
\item[1)] A map $U:L^{w,1}(\A)\ra L^{w,1}(\B)$ is a surjective isometry if and only if there exists a Jordan $*$-isomorphism $\Phi:\tma \ra \nmb$, a unitary $u\in \B$ and an $\alpha>0$ such that $U(x)=\alpha u\Phi(x)$ for all $x\in L^{w,1}(\A)$ and 
\begin{eqnarray}
\psi(\nu(\Phi(p)))=\frac{1}{\alpha} \psi(\tau(p)) \qquad \forall p\in \paf. \label{e0 T5}
\end{eqnarray}
\item[2)] If $\Phi:\A \ra \B$ is a Jordan $*$-isomorphism, $u\in \B$ is a unitary and there exists an $\alpha>0$ such that (\ref{e0 T5}) holds, then letting $U_0(x):=\alpha u\Phi(x)$ for all $x\in \A \cap L^{w,1}(\A)$ yields a map which can be extended to a surjective isometry $U:L^{w,1}(\A)\ra L^{w,1}(\B)$.
\end{itemize}
\begin{proof}
Suppose $U$ is a surjective isometry. We have seen that if $p\in \paf$, then $U(p)=\alpha_p v_{p}$ for some $v_{p}\in \mathcal{V}(\A)_f$.  It follows that $\nu(s(U(p)))=\nu(|v_{p}|)<\infty$ whenever $p\in \paf$. If $p,q\in \paf$ with $pq=0$, then $U(p)^*U(q)=0=U(p)U(q)^*$, by Proposition \ref{TC5.0 C1}. It follows by \cite[Theorem 5.3]{key-dJ19a}, that there exist a Jordan $*$-isomorphism $\Phi$, a unitary operator $u\in \B$ and a positive operator $b\in \nmb$ such that $U(x)=ub\Phi(x)$ for all $x\in \A \cap L^{w,1}(\A)$.

In order to show that this structural description can be improved, we start by showing that if $0\neq p,q\in \paf$ are arbitrary projections, then $\alpha_p=\alpha_q$. If $q<p$, then $0\neq p-q$ and $q(p-q)=0$. It follows by Proposition \ref{TL5.9 C1} that $\alpha_q=\alpha_{p-q}=\alpha_{q+(p-q)}=\alpha_p$. If $p \nleq q$ and $q\nleq p$, then let $m=p \vee q$. It follows that $p<m$, $q<m$ and $m\in \paf$. By what has been shown already this implies that $\alpha_p=\alpha_m=\alpha_q$. There therefore exists an $\alpha \in \mathbb{R}^+$ such that $U(p)=\alpha v_{p}$ holds for any $0\neq p\in \paf$.

Next, we show that $b=\alpha \id$. By considering \cite[Remark 5.4]{key-dJ19a}, we note that $\Phi(p)=s(U(p))$ for every $p\in \paf$ and the positive operator $b$ is constructed using the spectral projections of the positive operators $b_{p}:=|U(p)|$ ($p\in \paf$) by defining $e^b(\lambda,\infty):=\sotlimx{p\in \paf} e^{b_p}(\lambda,\infty)$ and constructing a positive operator from the resolution of the identity obtained in this way. In our present setting we have that \[b_{p}=|U(p)|=|\alpha v_{p}|=\alpha|v_{p}|=\alpha s(U(p))=\alpha \Phi(p),\] where we have used the fact that $\alpha_p=\alpha$ for all $p \in \paf$. This implies that if $b_{p}=\int_0^{\infty} \lambda de^{b_p}(\lambda)$ is the spectral representation of $b_{p}$, then \[e^{b_p}(\lambda, \infty)= \begin{cases} \Phi(p) & \text{if} \,\, 0\leq \lambda < \alpha \\
 0 & \text{if}\,\, \lambda \geq \alpha \end{cases} \] Furthermore, $\Phi$ is a Jordan $*$-isomorphism and hence normal and unital. It follows that $\sotlimx{p\in \paf} \Phi(p)=\Phi(\id)=\id$ and so \[e^b(\lambda,\infty):=\sotlimx{p\in \paf} e^{b_p}(\lambda,\infty)= \begin{cases} \id & \text{if} \,\, 0\leq \lambda < \alpha \\
 0 & \text{if}\,\, \lambda \geq \alpha \end{cases} \] It follows that $b=\int_0^\infty \lambda de^b(\lambda)=\alpha \id$ and hence $U(x)= \alpha u\Phi(x)$ for all $x\in \A \cap L^{w,1}(\A)$.  Furthermore, for any $p\in \paf$, we have 
\begin{eqnarray*}
&\psi(\tau(p))= \lwnorma{p} = \lwnormb{U(p)} =\lwnormb{\alpha u \Phi(p)} =\alpha \psi(\nu(\Phi(p))).&
\end{eqnarray*}

To obtain the desired global representation, we need to extend $\Phi$. In order to use \cite[Proposition 4.1]{key-Lab13} we need to show that $\nu\circ \Phi$ is $\epsilon-\delta$ absolutely continuous with respect to $\tau$. Let $\epsilon > 0$. Using the continuity of $\psi$ and $\psi^{-1}$ at $0$,  we can  find a $s_0>0$ such that $\psi^{-1}(s)<\epsilon$ if $0\leq s<s_0$, and a $\delta>0$ such that $\psi(t)<\alpha s_0$ if $t<\delta$. It follows that if $\tau(p)<\delta$, then $\psi(\nu\circ \Phi(p))=\frac{1}{\alpha}\psi(\tau(p))<s_0$ and therefore $\nu\circ \Phi(p)<\epsilon$. It follows that $\Phi$ can be extended to a continuous (with respect to the measure topologies) Jordan $*$-homomorphism from $\tma$ into $\nmb$. To show that this extension is bijective, we note that since $\Phi^{-1}$ is a Jordan $*$-isomorphism, it maps the projection lattice of $\B$ onto the projection lattice of $\A$. Using (\ref{e0 T5}), it is then easily checked that $\psi(\tau\circ \Phi^{-1}(q))=\alpha \psi(\nu(q))$ for every $q\in \pbf$. One can therefore similarly extend $\Phi^{-1}$ to a Jordan $*$-isomorphism from $\nmb$ into $\tma$. As in the proof of Lemma \ref{L2 17/09/19}, we can then show that the extension of $\Phi$ is, in fact, a Jordan $*$-isomorphism from $\tma$ onto $\nmb$. We will also denote this extension $\Phi$.

To show that the desired representation holds globally let $x\in L^{w,1}(\A)$. Since $L^{w,1}(\A)$ has order continuous norm, $L^{w,1}(\A)\cap \A$ is dense in $L^{w,1}(\A)$ and therefore there exists a sequence $(x_n)_{n=1}^\infty \subset L^{w,1}(\A)$ such that $x_n \rax{w,1} x$. It follows that $U(x_n)\rax{w,1} U(x)$ and therefore $U(x_n)\mtopc U(x)$, since $L^{w,1}(\B)$ is continuously embedded in $\nmb$. Since $x_n \rax{w,1} x$, we similarly have that $x_n \mtopc x$ and therefore $\Phi(x_n)\mtopc \Phi(x)$. Since $\nmb$ is a topological $*$-algebra with respect to the measure topology, it follows that $\alpha u\Phi(x_n) \mtopc \alpha u \Phi(x)$. However, $\alpha u\Phi(x_n)=U(x_n) \mtopc U(x)$ and so $U(x)=\alpha u \Phi(x)$ for every $x\in L^{w,1}(\A)$.

To prove 2) suppose that $\Phi:\A \ra \B$ is a Jordan $*$-isomorphism for which (\ref{e0 T5}) holds for some $\alpha>0$ and suppose $u\in \B$ is a unitary operator. Let $U_0(x):=\alpha u\Phi(x)$ for all $x\in \A\cap L^{w,1}(\A)$. We start by showing that $U_0$ is isometric on $\mathcal{G}(\A)_f$. Suppose $x=\sumx{i=1}{n} \alpha_i p_i \in \mathcal{G}(\A)_f$, with $|\alpha_1|>|\alpha_2|>...>|\alpha_n|$ and $p_i\in \paf$ for every $i$, with $p_i p_j=0$ if $i\neq j$. Then $|x|=\sumx{i=1}{n} |\alpha_i|p_i$. Furthermore, since $p_i p_j=0$ if $i\neq j$, it is easily checked that $x^*x =xx^* $. Note also that $|x|$ can be written in the form $|x|=\sumx{j=1}{n}\beta_j q_j$, where $q_1\leq q_2\leq ...\leq q_n$ and the $q_i$'s are projections. It is easily checked that $\svf{x}=\svf{|x|}=\sumx{j=1}{n}\beta_j \svf{q_j}$ and that $|\Phi(x)|=\Phi(|x|) =\sumx{j=1}{n}\beta_j \Phi(q_j)$. Since $\Phi$ is a Jordan $*$-isomorphism, $\Phi$ is positive and so $\Phi(q_1)\leq \Phi(q_2)\leq...\leq \Phi(q_n)$. Furthermore, these are all projections and so $\svf{\Phi(x)}=\svf{|\Phi(x)|}=\sumx{j=1}{n}\beta_j \chi_{[0,\nu(\Phi(q_j)))}$. Therefore 
\begin{eqnarray*}
&\lwnormb{U_0(x)}=\lwnormb{\alpha u\Phi(x)}  =\alpha\int_0^\infty \svft{\Phi(x)}{t}w(t)dt& \\
&=\alpha \sumx{j=1}{n}\beta_j \psi(\nu(\Phi(q_j))) = \alpha \sumx{j=1}{n} \frac{\beta_j}{\alpha} \psi(\tau(q_j))=\int_0^{\infty} \svft{x}{t}w(t) dt =\lwnorma{x}.&
\end{eqnarray*}
Even though $\mathcal{G}(\A)_f$ is dense in $L^{w,1}(\A)$ (see Remark \ref{RT2 29/01/15}), it need not be a subspace of $L^{w,1}(\A)$ and so care needs to be taken when extending $U_0$. Motivated by this, we show next that $U_0$ is isometric on $\mathcal{F}(\tau)$. It follows from the Spectral Theorem that if $x\in \mathcal{F}(\tau)^{sa}$, then there exists a sequence $\seq{x}$ in $\mathcal{G}(\A)_f^{sa}$ such that $x_n \rax{\A} x$ and $s(x_n)\leq s(x)$ for every $n\in \mathbb{N}^+$. Since $\Phi$ is a Jordan $*$-isomorphism and hence continuous, we have that $\Phi(x_n)\rax{\B} \Phi(x)$. If $y\in \A^{sa}$, then it is easily checked that $\Phi(s(y))=s(\Phi(y))$, since $y$ commutes with $s(y)$. It follows that $s(\Phi(x_n))\leq s(\Phi(x))$ for every $n$. Furthermore $\tau(s(x))<\infty$ implies that $\psi(\nu(\Phi(s(x))))=\tfrac{1}{\alpha}\psi(\tau(s(x)))<\infty$ and therefore $\nu(s(\Phi(x)))=\nu(\Phi(s(x)))<\infty$. A straightforward calculation then shows that $\Phi(x_n) \rax{L^{w,1}} \Phi(x)$ and similarly $x_n \rax{L^{w,1}} x$. Since $\lwnorm{\Phi(x_n)}=\lwnorm{x_n}$ for each $n$, we have that $\lwnorm{u\Phi(x)}=\lwnorm{\Phi(x)}=\lwnorm{x}$.

We have shown that $U_0$ is isometric on $\mathcal{F}(\tau)^{sa}$. It therefore has a unique isometric extension from $ L^{w,1}(\A)^{sa}$ into $L^{w,1}(\B)$, which can be further extended to a continuous map $U: L^{w,1}(\A)\ra L^{w,1}(\B)$, using the fact that every $x\in  L^{w,1}(\A)$ has a unique decomposition $x=x_1+ix_2$, with $x_1,x_2\in  L^{w,1}(\A)^{sa}$. This map $U$ is in fact an isometry, since it is isometric on the dense subset $\mathcal{G}(\A)_f$.

We show that $U$ is surjective. To do so, we will start by showing that $\Phi(\paf)=\pbf$. If $p\in \paf$, then $\Phi(p)\in \pb$. Furthermore, using (\ref{e0 T5}), we have that $\psi(\nu(\Phi(p)))=\tfrac1{\alpha} \psi(\tau(p)) <\infty$ and so $\Phi(\paf)\subseteq \pbf$. Suppose $q\in \pbf$. Since $\A$ is semi-finite and $p=\Phi^{-1}(q)$ is a projection, there exists a net $\{p_{\lambda}\}_{\lambda \in \Lambda}$ in $\paf$ such that $p_\lambda \uparrow p$. Therefore $\Phi(p_\lambda) \uparrow \Phi(p)$, since $\Phi$ is normal. This implies that $\psi(\nu(\Phi(p_\lambda))) \uparrow \psi(\nu(\Phi(p)))$ since $\nu$ is normal and $\psi$ is increasing and continuous. We can similarly show that $\frac{1}{\alpha} \psi(\tau(p_\lambda))\uparrow \frac{1}{\alpha} \psi(\tau(p))$, but  $\psi(\nu(\Phi(p_\lambda)))=\frac{1}{\alpha} \psi(\tau(p_\lambda))$ for all $\lambda$ and therefore $\frac{1}{\alpha} \psi(\tau(p))=\psi(\nu(\Phi(p)))=\psi(\nu(q))<\infty$. It follows that $p\in \paf$ and therefore $\pbf\subseteq \Phi(\paf)$. Since $\Phi(\paf)=\pbf$,  it follows that $\Phi(\mathcal{G}(\A)_f)=\mathcal{G}(\B)_f$. Suppose $y\in L^{w,1}(\B)$. Then $\tfrac{1}{\alpha} u^*y\in L^{w,1}(\B)$. Using the density of $\mathcal{G}(\B)_f$ in $L^{w,1}(\B)$ we can find a sequence $\seq{y}$ in $\mathcal{G}(\B)_f$ such that $y_n \rax{L^{w,1}(\B)} \tfrac{1}{\alpha} u^*y$, and for each $n$, let $x_n\in \mathcal{G}(\A)_f$ be such that $\Phi(x_n)=y_n$. Using the fact that $U$ is isometric on $\mathcal{F}(\tau)$ and $U(x_n)=\alpha u\Phi(x_n)$ for every $n$, one can show that $\seq{x}$ is Cauchy in $L^{w,1}(\A)$. Therefore $x_n \rax{L^{w,1}(\A)} x$, for some $x\in L^{w,1}(\A)$. It is easily checked that $U(x)=y$.  

It is easily checked that any Jordan $*$-isomorphism from $\tma$ onto $\nmb$ restricts to a Jordan $*$-isomorphism from $\A$ onto $\B$ and therefore the sufficiency of part  1) follows from 2).
\end{proof}
\end{thm}

In the following remark we show that in the finite setting Theorem \ref{T5 C1} reduces to \cite[Theorem 5.1]{key-Chilin89}.

\begin{rem} Suppose $(\A,\tau)$ and $(\B,\nu)$ are von Neumann algebras equipped with faithful normal finite traces such that $\tau(\id)=1=\nu(\id)$ and suppose $w:[0,1]\ra [0,\infty)$ is a strictly decreasing weight function with $\psi(1)=\int_0^1 w(t)dt=1$. We show that if $\Phi:\A \ra \B$ is a Jordan $*$-isomorphism, then condition (\ref{e0 T5}) is equivalent to $\Phi$ being trace-preserving. It is clear that if $\Phi$ is trace-preserving, then (\ref{e0 T5}) holds for $\alpha=1$. Suppose (\ref{e0 T5}) holds for some $\alpha>0$. In particular, $\alpha\psi(\nu(\Phi(\id)))=\psi(\tau(\id))$ and therefore $\alpha= \frac{\psi(\tau(\id))}{\psi(\nu(\id))}=1$, since $\Phi(\id)=\id$, $\tau(\id)=1=\nu(\id)$ and $\psi(1)=1$. This implies that $\psi(\nu(\Phi(p)))= \psi(\tau(p))$ for all $ p\in \pa$. Since $\psi$ is strictly increasing and hence injective, it follows that $\nu(\Phi(p))=\tau(p)$ for all  $p\in \pa$ and, therefore, $\Phi$ is trace-preserving, by Lemma \ref{L2 17/09/19}.
\end{rem}

In \cite[$\S6$]{key-Huang19} it is shown that the Lorentz spaces $\Lambda_w^p$ have strictly log monotone norm and hence that the structure of a positive surjective isometry on $\Lambda_w^p(\A)$ follows from \cite[Corollary 5.5]{key-Huang19}. Since $\Lambda_w^p(\A)=L^{w,1}(\A)$ if $p=1$, the result in \cite[$\S6$]{key-Huang19} holds for more general Lorentz spaces, but under the additional assumption that the isometry is positive. It is worth mentioning that even though \cite[Corollary 5.5]{key-Huang19} follows from \cite[Corollary 5.5]{key-Huang19}, which does not require positivity of the map, the proof of \cite[Corollary 5.5]{key-Huang19} also uses \cite[Proposition 4.8]{key-Huang19}, which demonstrates how the disjointness preserving property of a map may be obtained from the positivity of that map.

\section{Surjective isometries from $L^1+L^\infty(\A)$ onto  $L^1+L^\infty(\B)$} \label{S4}

In this section we show that surjective isometries from $L^1+L^\infty(\A)$ onto  $L^1+L^\infty(\B)$ can be characterized in terms of a unitary operator $u\in \B$ and a Jordan $*$-isomorphism $\tilde{\Phi}$ from $\tma$ onto $\nmb$. The key component in describing the structure of a surjective isometry $U:L^1+L^\infty(\A)\ra L^1+L^\infty(\B)$ is showing that $U\restriction_\A$ is an $L^\infty$-isometry from $\A$ onto $\B$, whose structure can therefore be described, using Kadison's Theorem (\cite[Theorem 7]{key-K51}). We start by proving a non-commutative analogue of \cite[Theorem 4]{key-Hudzik93}, which describes the extreme points of the unit ball of $L^1+L^\infty(\A)$. 

\begin{prop} \label{PL1 09/09/19}
Suppose $(\A,\tau)$ is a non-atomic semi-finite von Neumann algebra and let $E(0,\tau(\id))=L^1+L^{\infty}(0,\tau(\id))$. If $\tau(\id)\leq 1$, then $B_{E(\A)}$ does not have any extreme points. If $\tau(\id)>1$, then $x\in S_{E(\A)}$ is an extreme point of $B_{E(\A)}$ if and only if $x$ is a partial isometry with $\tau(|x|)=\tau(\id)$ and $n(x)\A n(x^*)=\{0\}$.
\begin{proof}
If $x\in S_{E(\A)}$ is an extreme point of $B_{E(\A)}$, then $\svf{x}\restriction_{(0,\tau(\id))}$ is an extreme point of $B_{E(0,\tau(\id))}$, by Theorem \ref{T2.2 Czer12}. Since $B_{E(0,\tau(\id))}$ does not have extreme points if $\tau(\id)\leq 1$ (see \cite[Theorem 4]{key-Hudzik93}), it follows that $B_{E(\A)}$ does not have extreme points if $\tau(\id)\leq 1$. Suppose $\tau(\id)>1$ and  $x\in S_{E(\A)}$ is an extreme point of $B_{E(\A)}$. Then $\svft{x}{t}=1$ for all $t\in (0,\tau(\id))$ using Theorem \ref{T2.2 Czer12} and \cite[Theorem 4]{key-Hudzik93}. Therefore $x=v$ for some partial isometry $v$ with $\tau(|v|)=\tau(\id)$ (see Lemma \ref{LP1 20/05/15}). If $\tau(\id)=\infty$, then $\limx{t\ra \infty} \svft{x}{t}=1\neq 0$ and so $n(x)\A n(x^*)=\{0\}$, by Theorem \ref{T2.2 Czer12}. If $1<\tau(\id)<\infty$, then $s(x)=|v|=\id$ (since $\tau(\id - |v|)=0$ and $\id - |v|\geq 0$) and so $n(x)\A n(x^*)=\{0\}$, since $n(x)=\id-s(x)=0$.

If $x$ is a partial isometry with $\tau(|x|)=\tau(\id)$ and $n(x)\A n(x^*)=\{0\}$, then $\svft{x}{t}=1$ for all $t\in (0,\tau(\id))$. It follows by \cite[Theorem 4]{key-Hudzik93} that $\svf{x}$ is an extreme point of $B_{E(0,\tau(\id))}$. If $1<\tau(\id)<\infty$, then $\limx{t\ra\infty}\svft{x}{t}=0$; if $\tau(\id)=\infty$, then $|x|=\id =\svft{x}\infty s(x)$ (and, by assumption,  $n(x)\A n(x^*)=\{0\}$), and so in both cases we are able to apply Theorem \ref{T2.2 Czer12} to conclude that $x$ is an extreme point of $B_{E(\A)}$.
\end{proof}
\end{prop} 

\begin{rem}\label{R1 09/09/19}
Since the set of extreme points of $B_{L^\infty(0,\alpha)}$ coincides with the set of extreme points of $B_{L^1+L^\infty(0,\alpha)}$ for $\alpha>1$ (see \cite[Theorem 4]{key-Hudzik93}), it follows, using the techniques employed above, that if $(\A,\tau)$ is a non-atomic semi-finite von Neumann algebra with $\tau(\id)>1$, then the set of extreme points of $B_{L^\infty(\A)}$ coincides with the set of extreme points of $B_{L^1+L^\infty(\A)}$.
\end{rem}

\begin{thm} \label{T1 18/09/19}
Suppose $(\A,\tau)$ and $(\B,\nu)$ are non-atomic semi-finite von Neumann algebras with $\tau(\id),\nu(\id)>1$. If $U:L^1+L^\infty(\A) \ra L^1+L^\infty(\B)$ is a surjective isometry, then there exist a unique trace-preserving Jordan $*$-isomorphism $\tilde{\Phi}$ from $\tma$ onto $\nmb$ and a unitary operator $u=U(\id)\in \B$ such that $U(x)=u\tilde{\Phi}(x)$ for all $x\in L^1+L^\infty(\A)$. Furthermore, $U\restriction_\A$ is an $L^\infty$-isometry from $\A$ onto $\B$ and $U\restriction_{L^1(\A)}$ is an $L^1$-isometry from $L^1(\A)$ onto $L^1(\B)$. Conversely, if $\tilde{\Phi}:\tma \ra \nmb$ is a trace-preserving Jordan $*$-isomorphism and $u\in \B$ is a unitary operator, then letting $U(x)=u\tilde{\Phi}(x)$ for $x\in L^1+L^\infty(\A)$ yields a surjective isometry from $L^1+L^\infty(\A)$ onto $L^1+L^\infty(\B)$.
\begin{proof}
Let $U:L^1+L^\infty(\A) \ra L^1+L^\infty(\B)$ be a surjective isometry. We start by showing that $U(\A)=\B$ and $\norminf{U(x)}=\norminf{x}$ for every $x\in \A$. Suppose $x\in B_\A$. Then $x=\tfrac{1}{2}(v+w)$ for some $v,w\in \text{ext}\,(B_\A)$, by \cite[Theorem 3]{key-Choda70}. Since $\text{ext}\,(B_{L^\infty(\A)})=\text{ext}\,(B_{L^1+L^\infty(\A)})$ (see Remark \ref{R1 09/09/19}) and $U$ preserves extreme points, we have that $U(v),U(w)\in \text{ext}\,(B_{L^1+L^\infty(\B)})=\text{ext}\,(B_{L^\infty(\B)})$ and hence $U(x)=\tfrac{1}{2}\left(U(v)+U(w)\right)\in B_{L^\infty(\B)}$. Since $U^{-1}$ is also an isometry, we have that $U(B_\A)=B_\B$ and hence $U(\A)=\B$. Suppose $0\neq x \in \A$. Since $U(B_\A)=B_\B$, we have that $\norminf{U\left(\frac{x}{\norminf{x}}\right)}\leq 1$ and therefore $\norminf{U(x)}\leq \norminf{x}$. Since we also have that $U^{-1}(B_\B)=B_\A$, we obtain $\norminf{x}=\norminf{U^{-1}(U(x))}\leq \norminf{U(x)}$.

It follows from \cite[Theorem 7]{key-K51} that there exists a Jordan $*$-isomorphism $\Phi$ from $\A$ onto $\B$ and a unitary operator $u=U(\id)$ such that $U(x)=u\Phi(x)$ for every $x\in \A$. To show that $\Phi$ can be extended we demonstrate that $\Phi$ is trace-preserving on projections with finite trace. Suppose $p\in \paf$ with $\tau(p)<1$. Then using the fact that $\norm{x}_{1+\infty}=\int_0^1 \svft{x}{t}dt$ (see \cite[p.230]{key-Dodds14}), we obtain
\begin{eqnarray*}
\tau(p)&=&\min\{1,\tau(p)\}=\pltnorm{p}=\plnnorm{U(p)}\\
&=&\plnnorm{u\Phi(p)}=\plnnorm{\Phi(p)}=\min\{1,\nu(\Phi(p))\}.
\end{eqnarray*}
Therefore $\nu(\Phi(p))=\tau(p)$. If $p\in \paf$ with $\tau(p)\geq 1$, then using the non-atomicity of $\A$ we can find $(p_n)_{n=1}^k\subseteq \paf$ such that $p_np_m=0$ if $n\neq m$, $\sumx{n=1}{k}p_n=p$ and $\tau(p_n)<1$ for every $n$. This enables one to conclude that $\nu(\Phi(p))=\tau(p)$. It now follows from Lemma \ref{L2 17/09/19} that $\Phi$ has a unique extension to a Jordan $*$-isomorphism $\tilde{\Phi}$ from $\tma$ onto $\nmb$ and $\tilde{\Phi}$ is trace-preserving on $\tma$. We show that if $x\in L^1+L^\infty(\B)$, then $U(x)=u\tilde{\Phi}(x)$. Since $L^1\cap L^\infty(\A)$ is dense in $L^1(\A)$ and $L^1(\A)$ is continuously embedded in $L^1+L^\infty(\A)$, we have that $L^1(\A)\subseteq \overline{L^1\cap L^\infty(\A)}^{L^1+L^\infty}$. This implies that $\A$ is dense in $L^1+L^\infty(\A)$ with respect to the $L^1+L^\infty$-norm. We can therefore find a sequence $\seq{x}\subseteq \A$ such that $x_n \plc x$ and hence $U(x_n) \plc U(x)$. Using the continuity of the embedding of an $L^1+L^\infty$-space into the corresponding space of trace-measurable operators (see \cite[Proposition 20]{key-Dodds14}), it follows that $x_n \mtopc x$ and $U(x_n) \mtopc U(x)$. Therefore $\tilde{\Phi}(x_n) \mtopc \tilde{\Phi}(x)$, by \cite[Theorem 3.10]{key-Weigt09} and hence $u\tilde{\Phi}(x_n) \mtopc u\tilde{\Phi}(x)$. However, $U(x_n)=u\tilde{\Phi}(x_n)$ for every $n$ and therefore $U(x)=u\tilde{\Phi}(x)$.

To demonstrate the uniqueness of $u$ and $\tilde{\Phi}$, suppose $v\in \mathcal{B}$ is a unitary and $\tilde{\Psi}$ is a Jordan $*$-isomorphism from $\tma$ onto $\nmb$ such that $U(x)=v\tilde{\Psi}(x)$ for every $x\in L^1+L^\infty(\A)$. It follows by \cite[Lemma 4.2]{key-Weigt09} that $\tilde{\Psi}(\A)\subseteq \B$ and $\tilde{\Psi}^{-1}(\B)\subseteq \A$. Therefore $\tilde{\Psi}\restriction_{\A}$ is a Jordan $*$-isomorphism from $\A$ onto $\B$ and hence unital. Since $\tilde{\Phi}$ is likewise unital, we have that $u=U(\id)\tilde{\Phi}(\id)=U(\id)=v\tilde{\Psi}(\id)=v$ and furthermore, for any $x\in L^1+L^\infty(\A)$, \[\tilde{\Psi}(x)=u^*u\tilde{\Psi}(x)=u^*U(x)=u^*u\tilde{\Phi}(x)=\tilde{\Phi}(x).\] 
Since $\A$ is dense in $\tma$ with respect to the measure topology, $\tilde{\Phi}$ and $\tilde{\Psi}$ agree on $\A$ and both are continuous with respect to the measure topology (see \cite[Theorem 3.10]{key-Weigt09}), we have that $\tilde{\Phi}=\tilde{\Psi}$. 
Finally, the fact that $U\restriction_{L^1(\A)}$ is an $L^1$-isometry from $L^1(\A)$ onto $L^1(\B)$ follows from \cite[Theorem 2 and its Corollary on p.49]{key-Y81} (see also \cite[Remark 2.2.5]{key-dJ17}).

Conversely, suppose $\tilde{\Phi}:\tma \ra \nmb$ is a trace-preserving Jordan $*$-isomorphism and $u\in \B$ is a unitary operator. Let $U(x)=u\tilde{\Phi}(x)$ for $x\in L^1+L^\infty(\A)$. As before, we have that $\tilde{\Phi}(\A)=\B$ and that $\tilde{\Phi}\restriction_\A$ is a Jordan $*$-isomorphism from $\A$ onto $\B$. Therefore $\Phi\restriction_{\A}$ is $L^\infty$-isometric. It follows from the fact that $\tilde{\Phi}$ (and hence also $\tilde{\Phi}^{-1}$) is trace-preserving that $\tilde{\Phi}\restriction_{L^1(\A)}$ is an $L^1$-isometry from $L^1(\A)$ onto $L^1(\B)$. It is then easily checked that $\tilde{\Phi}(L^1+L^\infty(\A))=L^1+L^\infty(\B)$ and $\plnnorm{u\tilde{\Phi}(x)}=\plnnorm{\tilde{\Phi}(x)}=\pltnorm{x}$ for every $x\in L^1+L^\infty(\B)$. 
\end{proof}
\end{thm}

\begin{rem}
Suppose $(\A,\tau)$ and $(\B,\nu)$ are semi-finite von Neumann algebras. If $\tau(\id_\A)<\infty$, then $L^1+L^\infty(\A)=L^1(\A)$. If, in addition $\tau(\id_\A)\leq 1$, then it is easily checked that $\norm{x}_{1+\infty}=\norm{x}_{1}$ for every $x\in L^1+L^\infty(\A)$ (see \cite[p.145]{key-Hudzik93} for the proof of the corresponding claim in the commutative setting). Consequently, if $\tau(\id_{\A}),\nu(\id_{\B})\leq 1$, then surjective isometries from $L^1+L^\infty(\A)$ onto $L^1+L^\infty(\B)$ are characterized by Yeadon's Theorem (\cite[Theorem 2]{key-Y81}). 
\end{rem}

If $U:L^1+L^\infty(\A) \ra L^1+L^\infty(\B)$ is a positive surjective isometry, then $U(\id)=\id$, since $U(\id)\geq 0$ and $U(\id)$ is unitary by the previous theorem. We therefore obtain the following characterization of positive surjective isometries from $L^1+L^\infty(\A)$ onto $L^1+L^\infty(\B)$.

\begin{cor}
Suppose $(\A,\tau)$ and $(\B,\nu)$ are non-atomic semi-finite von Neumann algebras with $\tau(\id),\nu(\id)>1$. Then $U:L^1+L^\infty(\A) \ra L^1+L^\infty(\B)$ is a positive surjective isometry if and only if $U$ is the restriction of a trace-preserving Jordan $*$-isomorphism $\tilde{\Phi}$ from $\tma$ onto $\nmb$. 
\end{cor}

\section{Positive surjective isometries from $L^1\cap L^\infty (\A)$ onto $L^1 \cap L^\infty (\B)$ } \label{Intersection isometries}

In this section we show that positive surjective isometries from  $L^1\cap L^\infty(\A)$ onto $L^1\cap L^\infty(\B)$ can be characterized as the restrictions of trace-preserving Jordan $*$-isomorphisms.  Since we are not assuming that $\A=\B$ nor that the identities on these spaces necessarily have infinite trace, the main result of this section yields new information even in the commutative setting (see \cite[Theorem 1]{key-Grz92b}). The following characterization of the extreme points of the unit ball of $L^1\cap L^\infty(\A)$ will play a significant role in obtaining the desired structural description and can be proved using \cite[Corollary 1]{key-Grz92} and a similar technique (with a few minor adjustments) to the one employed in the proof of Proposition \ref{PL1 09/09/19}.

\begin{prop} \label{P2 10/09/16}
Suppose $(\A,\tau)$ is a non-atomic semi-finite von Neumann algebra and let $E(0,\tau(\id))=L^1\cap L^\infty(0,\tau(\id))$. Then $x$ is an extreme point of $B_{E(\A)}$ if and only if $x=v$ for some partial isometry $v\in \A$ with $\tau(|v|)=\text{min}\,\{1,\tau(\id)\}$.
\end{prop}

For the remainder of this section we will assume that $(\A,\tau)$ and $(\B,\nu)$ are non-atomic semi-finite von Neumann algebras with $\tau(\id),\nu(\id)\geq 2$. Suppose $U$ is a positive surjective isometry from $L^1\cap L^\infty(\A)$ onto $L^1 \cap L^\infty (\B)$. To show that $U$  has the desired structure we will use the fact that any linear positive normal map from $\mathcal{F}(\tau)$ into $\B$ which is square-preserving on self-adjoint elements can be extended to a normal Jordan $*$-homomorphism from $\A$ into $\B$ (see \cite[Theorem 4.5]{key-dJ18b}). It follows from the positivity and surjectivity of $U$ that $U$ is normal (see \cite[Lemma 3.1]{key-dJ19a}) and therefore we need to show that $U$ is square-preserving on self-adjoint elements in $\mathcal{F}(\tau)$. We start by showing that $U$ maps projections with finite trace onto projections and that $U$ is projection disjointness preserving.

\begin{lem} \label{L1 14/09/16}
If $p\in \paf$, then $U(p)\in \pb$. Furthermore, if $p,q\in \paf$ with $pq=0$, then $U(p)U(q)=0$.
\begin{proof}
Suppose $p\in \paf$ with $\tau(p)=1$. By Proposition \ref{P2 10/09/16}, $p$ is an extreme point of the unit ball of $L^1\cap L^\infty (\A)$. Since $U$ is a surjective isometry, we have that $U(p)$ is an extreme point of the unit ball of $L^1 \cap L^\infty(\B)$ and hence $U(p)=v_{p}$ for some partial isometry $v_{p}\in \B$ with $\nu(|v_{p}|)=1$. Using the positivity of $U$ we have that $U(p)=v_p=|v_p|$ is a projection and  $\nu(U(p))=1$. 
 
Next, we show that $U(p)U(q)=0$ if $p,q\in \paf$ with $pq=0$ and $\tau(p)=1=\tau(q)$.  It is easily checked that $|p-q|=|p+q|$. Since $\intnorm{|y|}=\intnorm{y}$ for every $y\in L^1 \cap L^\infty (\A)$, it follows that 
\[ \intnorm{p-q}=\intnorm{p+q} =\text{max}\,\Bigl\{\norm{p+q}_{1},\norm{p+q}_\infty\Bigr\} =2.  \]
Therefore $\intnorm{U(p)-U(q)}=2$. Since $U(p)$ and $U(q)$ are projections, we have that $-\id \leq U(p)-U(q) \leq \id$ and $-U(q)\leq U(p)-U(q)\leq U(p)$. Therefore $\norm{U(p)-U(q)}_\infty \leq \norm{\id}_\infty= 1$  (see \cite[Proposition 4.2.8]{key-K1}) and $\norm{U(p)-U(q)}_{1}\leq \norm{U(p)+U(q)}_{1}$ (see \cite[Corollary 4]{key-Bik12}). Hence
\begin{eqnarray*}
2&=& \text{max}\,\Bigl\{\norm{U(p)-U(q)}_{1},\norm{U(p)-U(q)}_\infty \Bigr\} = \norm{U(p)-U(q)}_{1} \\
&\leq& \norm{U(p)+U(q)}_{1} \leq \norm{U(p)}_{1}+\norm{U(q)}_{1} =2, 
\end{eqnarray*}
since $U(p)$ and $U(q)$ are projections with $\nu(U(p))=1=\nu(U(q))$. It follows that 
\begin{eqnarray}
\norm{U(p)+U(q)}_{1}+\norm{U(p)-U(q)}_{1}=2\Bigl(\norm{U(p)}_{1}+\norm{U(q)}_{1}\Bigr),  \label{CE}
\end{eqnarray}
i.e. we have equality in Clarkson's inequality. Application of \cite[Theorem 11.4.3]{key-F2} yields $U(p)U(q)=0$. 

Next, suppose $p,q\in \paf$ with $pq=0$ and $0<\tau(p),\tau(q)<1$. Since $(\A,\tau)$ is non-atomic and $\tau(\id)>1$, there exist  $p_1,q_1\in \paf$ such that $p+p_1,q+q_1\in \paf$, $\tau(p+p_1)=1=\tau(q+q_1)$ and $(p+p_1)(q+q_1)=0$. It follows by what has been shown already that $U(p+p_1)$ and $U(q+q_1)$ are orthogonal projections. Furthermore, $0\leq U(p)\leq U(p+p_1)$ and therefore $s(U(p))\leq s(U(p+p_1))=U(p+p_1)$. Similarly, $r(U(q))\leq U(q+q_1)$ and thus $U(p)U(p+p_1)U(q+q_1)U(q)=U(p)U(q)$. However, $U(p+p_1)U(q+q_1)=0$ and so $U(p)U(q)=0$. Since $pp_1=0$ and $0<\tau(p),\tau(q)<1$, we therefore also have that $U(p)U(p_1)=0$. It is then easily checked (see \cite[Proposition B.1.32]{key-dJ17}, for example) that $s(U(p))s(U(p_1))=0$ and
\[s(U(p))+s(U(p_1))=s(U(p)+U(p_1))=s(U(p+p_1))=U(p+p_1):=e.\]
This implies that $s(U(p))\leq e$ and so for $\eta \in s(U(p))(H)$ we therefore have $\eta=e\eta=(U(p)+U(p_1))\eta=U(p)\eta$, since $U(p_1)\eta=U(p_1)s(U(p_1))s(U(p))\eta=0$. It follows that $U(p)$ is a projection and $U(p)=s(U(p))$.

Finally, for general $p,q\in \paf$ with $pq=0$, we use the non-atomicity of $(\A,\tau)$ to find  $(p_i)_{i=1}^k,(q_j)_{j=1}^n \subseteq \paf$ such that $\sumx{i=1}{k}p_i=p$ and $\sumx{j=1}{n}q_j=q$, $p_ip_j=0=q_iq_j$ if $i\neq j$ and $\tau(p_i),\tau(q_j)<1$ for all $i,j$. Since $p_i\leq p$ for each $i$ and $q_j\leq q$ for each $j$, we have that $p_iq_j=0$ for all $i,j$. Therefore $U(p_i)U(q_j)=0$ for all $i,j$ and so \[U(p)U(q)=\sumx{i=1}{k}U(p_i) \sumx{j=1}{n}U(q_j)=\sumx{i=1}{k}\sumx{j=1}{n}U(p_i)U(q_j)=0.\]
Furthermore, by what has been shown already, we have that $U(p_i)$ is a projection for each $i$ and these projections are mutually orthogonal. It follows that $U(p)=\sumx{i=1}{k}U(p_i)$ is a sum of mutually orthogonal projections and is therefore a projection (see \cite[Exercise 2.3.4]{key-Con1}).  
\end{proof}
\end{lem} 

\begin{cor}\label{C1 14/09/16}
If $x\in \mathcal{F}(\tau)^{sa}$, then $U(x^2)=U(x)^2$.
\begin{proof}
If $p\in \paf$, then $U(p^2)=U(p)=U(p)^2$, since $U(p)$ is a projection. If $x \in \mathcal{G}_f^{sa}$ (i.e. $x=\sumx{1=i}{n}\alpha_i p_i$ with $\alpha_i \in \mathbb{R}$, $p_i\in \paf$ for each $i$ and $p_i p_j=0$ if $i\neq j$), then $x^2=\sumx{1=i}{n}\alpha_i^2 p_i$ and therefore $U(x^2)=\sumx{1=i}{n}\alpha_i^2 U(p_i)= \left(\sumx{1=i}{n}\alpha_i U(p_i)\right)^2=U(x)^2$, using Lemma \ref{L1 14/09/16}. Suppose $x\in \mathcal{F}(\tau)^{sa}$. It follows from the Spectral Theorem that there exists $\seq{x}\subseteq \mathcal{G}^{sa}_f$ such that $x_n \rax{\A} x$ and $s(x_n)\leq s(x)$ for all $n\in \mathbb{N}^+$. Note that $x_n\raxl{L^{1}\cap L^\infty} x$ (see \cite[Proposition 2.1]{key-dJ19a}) and therefore $U(x_n)\raxl{L^1\cap L^\infty} U(x)$.  Furthermore, $\norm{U(x)-U(x_n)}_\B \leq \intnnorm{U(x_n)-U(x)}\ra 0$ and so $U(x_n)\rax{\B} U(x)$ (and so also $U(x_n)^2\rax{\B} U(x)^2$). Similarly, $U(x_n^2) \rax{\B} U(x^2)$. However, $U(x_n^2)=U(x_n)^2$ for every $n$ and so $U(x^2)=U(x)^2$. 
\end{proof}
\end{cor}

Equipped with these preliminary results we are now in a position to characterize positive surjective isometries from $L^1\cap L^\infty(\A)$ onto $L^1 \cap L^\infty (\B)$.

\begin{thm}\label{T 14/09/16}
Suppose $(\A,\tau)$ and $(\B,\nu)$ are non-atomic semi-finite von Neumann algebras with $2 \leq \tau(\id),\nu(\id)\leq \infty$. If $U:L^1\cap L^\infty(\A) \ra L^1 \cap L^\infty (\B)$ is a positive surjective isometry, then $U$ is the restriction to $L^1\cap L^\infty(\A)$  of a trace-preserving Jordan $*$-isomorphism from $\A$ onto $\B$. Conversely, if $U$  is a trace-preserving Jordan $*$-isomorphism from $\A$ onto $\B$, then $U$ is positive and maps $L^1\cap L^\infty(\A)$ isometrically onto $L^1\cap L^\infty (\B)$.
\begin{proof}
Suppose $U$ is a positive surjective isometry from $L^1\cap L^\infty(\A)$ onto $L^1 \cap L^\infty (\B)$. Since $U\restriction_{\mathcal{F}(\tau)}$ is linear, positive, normal and square-preserving on self-adjoint elements, it has a unique extension to a normal Jordan $*$-homomorphism $\Phi$ from $\A$ into $\B$, by \cite[Theorem 4.5]{key-dJ18b}. By \cite[Remark 4.6]{key-dJ18b} the injectivity of $U$ implies that $\norm{\Phi(x)}_\B=\norm{x}_\A$ for all $x\in \A^{sa}$. Furthermore, if $p\in \paf$ with $\tau(p)> 1$. Then $\tau(p)=\norm{p}_{1}=\intnorm{p} =\intnorm{U(p)}=\norm{U(p)}_{1}=\nu(U(p))$, since $\norm{p}_\infty=1 < \norm{p}_{1}$ and $\norm{U(p)}_\infty=1 <\intnorm{U(p)}$. If $p\in \paf$ with $\tau(p)\leq 1$, then there exists $p_1\in \pa$ with $p_1\leq p^\perp$ and $1<\tau(p_1)<\infty$. By what has been shown already, we have that $\nu(U(p)) =\nu(U(p+p_1))-\nu(U(p_1)) =\tau(p+p_1)-\tau(p_1) =\tau(p)$. It follows by Lemma \ref{L2 17/09/19} that $\Phi$ is trace-preserving.

If we can show that $\Phi$ is unital and $\Phi(p)\B \Phi(p) \subseteq \Phi(\A)$ for every $p\in \paf$, then it will follow from \cite[Proposition 6.2]{key-dJ18b} that $\Phi$ is a Jordan $*$-isomorphism from $\A$ onto $\B$. Since $(\B, \nu)$ is semi-finite, there exists an increasing net $\net{q}\subseteq \pbf$ such that $q_\lambda \uparrow \id$. Furthermore, $\pbf \subseteq L^1\cap L^\infty (\nu)$ and so  $U^{-1}(q_\lambda)$ is defined for each $\lambda$. Since $U^{-1}$ is a positive (see \cite[Lemma 3.1]{key-dJ19a}) surjective isometry, $\{U^{-1}(q_\lambda)\}_{\lambda \in \Lambda}$ is an increasing net of projections, by applying Lemma \ref{L1 14/09/16} to $U^{-1}$. It follows that $p_\lambda:=U^{-1}(q_\lambda)\uparrow p$ for some $p\in \pa$. Furthermore, $\Phi$ is normal and so $\Phi(p_\lambda)\uparrow \Phi(p)$, but \[\Phi(p_\lambda)=U(U^{-1}(q_\lambda))=q_\lambda \uparrow \id\] and so $\Phi(p)=\id$. Since $p\leq \id$, we have that $\Phi(p)\leq \Phi(\id)$. Since $\Phi$ maps projections onto projections (see \cite[Exercise 10.5.22(5)]{key-K2}), we have that $\Phi(\id)\leq \id=\Phi(p)$. It follows that $\Phi(\id)=\Phi(p)=\id$. Next, suppose $p\in \paf$. Then $\Phi(p)\in \pbf$, since $\Phi$ is trace-preserving. It follows that \[\Phi(p)\B \Phi(p) \subseteq \mathcal{F}(\nu) \subseteq L^1 \cap L^\infty(\nu)=U(L^1\cap L^\infty(\tau))\subseteq \Phi(\A).\] 

Conversely, suppose $U$ is a trace-preserving Jordan $*$-isomorphism from $\A$ onto $\B$. As in the proof of the  sufficiency part of Yeadon's Theorem (\cite[Theorem 2]{key-Y81}) one can show that $\norm{U(x)}_{1}=\norm{x}_{1}$ for all $x\in \A \cap L^1(\A)$. Furthermore, $U$ is an $L^\infty$-isometry and so for $x\in L^1 \cap L^\infty(\A)$ we have
\[\intnorm{U(x)}=\max\{\norm{U(x)}_1,\norm{U(x)}_\infty\}=\max\{\norm{x}_1,\norm{x}_\infty\}=\intnorm{x}.\] 
\end{proof}
\end{thm}

\begin{rem}
The structure of a positive surjective isometry $U:L^1\cap L^\infty(\A) \ra L^1 \cap L^\infty (\B)$ can also be described using other approaches. Having shown that such an isometry is projection disjointness preserving, it follows from \cite[Theorem 3.6]{key-Huang19} (see also \cite[Remark 4.12]{key-dJ19a}) that there exist a positive operator $b$ affiliated with $\B$ and a normal Jordan $*$-monomorphism $J$ from $\A$ onto a weakly closed $*$-subalgebra of $\B$ such that $U(x)=b J(x)$ for every $x\in L^1\cap L^\infty (\A)$. Alternatively, if one also shows that $U$ is finiteness preserving (using the non-atomicity of $\A$, the positivity of $U$ and the fact that $\nu(s(U(p))=\nu(U(p))=1$ for every projection $p\in \paf$ with $\tau(p)=1$), then one could apply \cite[Theorem 4.11]{key-dJ19a} to show that there exist a positive operator $a\in \nmb$ and a Jordan $*$-isomorphism $\Phi$ from $\A$ onto $\B$ such that $U(x)=a \Phi(x)$ for every $x\in L^1\cap L^\infty (\A)$. In both cases it would then be desirable to refine the initial representation.
\end{rem}

\begin{rem}
We make a few brief comments regarding the assumption $2 \leq  \tau(\id),\nu(\id)\leq \infty$. Note that if $\tau(\id_\A)<\infty$, then $L^1\cap L^\infty(\A)= \A$. If, in addition, $\tau(\id_\A)\leq 1$, then for any $x\in L^1\cap L^\infty(\A)= \A$, we have that $\norm{x}_1=\norm{x\id}_1\leq \norm{x}_\infty \tau(\id)\leq \norm{x}_\infty$ and therefore $\norm{x}_{1\cap \infty}=\max\{\norm{x}_1,\norm{x}_\infty\}=\norm{x}_\infty$. Consequently, if $\tau(\id_\A),\nu(\id_\B)\leq 1$, then surjective isometries from $L^1\cap L^\infty(\A)$ onto $L^1\cap L^\infty(\A)$  are described by Kadison's Theorem (\cite[Theorem 7]{key-K51}).  To show that $U$ is projection disjointness preserving in the proof of Lemma \ref{L1 14/09/16} we use a technique which relies on being able to use two orthogonal projections $p$ and $q$ for which we know the precise values of $\norm{U(p)}_1$ and $\norm{U(q)}_1$ in order to obtain equality in Clarkson's inequality (see (\ref{CE})). For $\tau(p)=1=\tau(q)$ we have access to the characterization of the extreme points of the unit ball, but not otherwise. For this reason we have assumed that $\tau(\id_\A)\geq \tau(p+q)\geq 2$. This also plays a role in showing that $U$ maps projections to projections and since we later apply this to $U^{-1}$, we also assume that $\nu(\id_\B)\geq 2$.  It is therefore an open problem to describe positive surjective isometries from $L^1\cap L^\infty(\A)$ onto $L^1 \cap L^\infty (\B)$ if $1<\tau(\id_\A)< 2$ or $1<\nu(\id_\B) < 2$.  
\end{rem}

\section*{Acknowledgments}

A significant part of this research was conducted during the first author's doctoral studies at the University of Cape Town and during his time at North West University. The first author would like to thank his Ph.D. supervisor, Dr Robert Martin, for his input and guidance and the NRF for funding towards this project in the form of scarce skills and grantholder-linked bursaries. Furthermore, the support of the DST-NRF Centre of Excellence in Mathematical and Statistical Sciences (CoE-MaSS) towards this research is hereby acknowledged. Opinions expressed and conclusions arrived at, are those of the authors and are not necessarily attributed to the CoE.

\end{document}